\documentclass[10pt,a4paper]{article}

\usepackage{amsmath, amsthm, amssymb, amsfonts}
\usepackage[T1]{fontenc}
\usepackage{graphicx}

\numberwithin{equation}{section}
\newtheorem{theorem}{Theorem}[section]
\newtheorem{lemma}[theorem]{Lemma}
\newtheorem{definition}[theorem]{Definition}
\newtheorem{corollary}[theorem]{Corollary}
\newtheorem{conjecture}[theorem]{Conjecture}

\newtheorem{remark}[theorem]{Remark}
\newtheorem{example}[theorem]{Example}

\DeclareSymbolFont{AMSb}{U}{msb}{m}{n}
\DeclareMathSymbol{\N}{\mathalpha}{AMSb}{"4E}
\DeclareMathSymbol{\R}{\mathalpha}{AMSb}{"52}
\DeclareMathSymbol{\Z}{\mathalpha}{AMSb}{"5A}
\DeclareMathSymbol{\D}{\mathalpha}{AMSb}{"44}
\DeclareMathSymbol{\X}{\mathalpha}{AMSb}{"58}
\DeclareMathSymbol{\s}{\mathalpha}{AMSb}{"53}
\newcommand{\M}{\mathsf{M}}
\newcommand{\m}{\mathsf{m}}

\newcommand{\V}{\mathsf{vol}}
\newcommand{\sk}{{\mathcal S}_\kappa}
\newcommand{\ck}{{\mathcal C}_\kappa}

\begin{document}

\title{Ricci Bounds for Euclidean and Spherical Cones}

\author{Kathrin Bacher, Karl-Theodor Sturm}

\date{}

\maketitle

\begin{abstract}
We prove generalized lower Ricci bounds for Euclidean and spherical cones over complete Riemannian manifolds. These cones are regarded as complete metric measure spaces. In general, they will be neither manifolds nor Alexandrov spaces. We show that the Euclidean cone over an $n$-dimensional Riemannian manifold whose Ricci curvature is bounded from below by $n-1$ satisfies the curvature-dimension condition $\mathsf{CD}(0,n+1)$ and that the spherical cone over the same manifold fulfills the curvature-dimension condition $\mathsf{CD}(n,n+1)$.

More generally, for each $N>1$ we prove that the condition $\mathsf{CD}(N-1,N)$ for a weighted Riemannian space is equivalent to the
condition $\mathsf{CD}(0,N+1)$ for its $N$-Euclidean cone as well as to the
condition $\mathsf{CD}(N,N+1)$ for its $N$-spherical cone.
\end{abstract}

\section{Introduction}

In two similar but independent approaches, the second author \cite{sa,sb} and Lott \& Villani \cite{lva,lv} presented a concept of generalized lower Ricci curvature bounds for metric measure spaces $(\M,\mathsf{d},\m)$.
The full strength of this concept appears if the condition $\mathsf{Ric}(\M,\mathsf{d},\m)\ge K$ is combined with a kind of upper bound $N$ on the dimension.
This leads to the so-called curvature-dimension condition $\mathsf{CD}(K,N)$ which can be formulated in terms of optimal transportation for each pair of numbers $K\in\mathbb{R}$ and $N\in[1,\infty)$.

A complete Riemannian manifold satisfies $\mathsf{CD}(K,N)$ if and only if its Ricci curvature is bounded from below by $K$ and its dimension from above by $N$.

A broad variety of geometric and functional analytic results can be deduced from the curvature-dimension condition $\mathsf{CD}(K,N)$. Among them are the Brunn-Minkowski inequality and the theorems by Bishop-Gromov, Bonnet-Myers and Lichnerowicz. Moreover, the condition $\mathsf{CD}(K,N)$ is stable under convergence with respect to the $\mathsf{L}^2$-transportation distance $\D$.

\subsection{Statement of the Main Results}

Let $\M$ be a complete $n$-dimensional Riemannian manifold (with Riemannian distance $\mathsf{d}$ and Riemannian volume $d\m=d\mathsf{vol}$). The \textit{Euclidean cone} $\mathsf{Con}(\M)=\M\,\, {}_r\!\!\times [0,\infty)$ over $\M$ is defined as the quotient of the product $\M\times [0,\infty)$ obtained by identifying all points in the fiber $\M\times\{0\}$. This point is called the origin $\mathsf{O}$ of the cone. It is equipped with a metric $\mathsf{d_{Con}}$ defined by the cosine formula
$$\mathsf{d_{Con}}((x,s),(y,t))=\sqrt{s^2+t^2-2st\cos(\mathsf{d}(x,y) \wedge \pi)},$$
and with a measure $m_n$ defined as the product $d \m_n(x,s):=d\m(x)\otimes s^nds$.

\begin{theorem} \label{Cheeger}
The Ricci curvature of $\M$ is bounded from below by $n-1$ and $\mathsf{diam}(M)\le\pi$ if and only if the metric measure space $(\mathsf{Con}(\M),\mathsf{d_{Con}},\m_n)$ satisfies the curvature-dimension condition $\mathsf{CD}(0,n+1)$.
\end{theorem}

Note that in dimensions $n\not=1$ the diameter bound $\mathsf{diam}(M)\le\pi$ is redundant: it follows from the Ricci bound.

\medskip

The heuristic interpretation of the assertion in the theorem is that the Euclidean cone -- regarded as a metric measure space -- has non-negative Ricci curvature in a generalized sense. Note that already in 1982, Cheeger and Taylor \cite{cha,chb} observed that the punctured Euclidean cone $\mathsf{Con}(\M)\setminus\{\mathsf{O}\}$ constructed over a compact $n$-dimensional Riemannian manifold $\M$ with $\mathsf{Ric}\geq n-1$ is an $(n+1)$-dimensional Riemannian manifold with $\mathsf{Ric}\geq 0$.
Note, however, that the sectional curvature might be unbounded from below (and above). Thus in general $\mathsf{Con}(\M)$ will \emph{not} be an Alexandrov space. Moreover,  $\mathsf{Con}(\M)$ in general is not a manifold and, of course, $\mathsf{Con}(\M)\setminus\{\mathsf{O}\}$ is not complete. In particular, the Ricci curvature in the classical sense is not defined in its singularity $\mathsf{O}$.

\medskip

Actually, we will prove a significantly more general result:
\begin{theorem} For any real number $N>1$, the $\mathsf{CD}(N-1,N)$ condition for a weighted Riemannian manifold is equivalent to the $\mathsf{CD}(0,N+1)$ condition for the associated $N$-Euclidean cone.
\end{theorem}

It is an open question whether analogous assertions hold true with an arbitrary metric measure space
$(\M,\mathsf{d},\m)$ in the place of the weighted Riemannian manifold $M$.
A partial result towards this conjecture  was derived by S. Ohta \cite{ob} for metric measure spaces satisfying the so-called \emph{measure contraction property} $\mathsf{MCP}(K,N)$, a property being slightly weaker than the curvature-dimension condition $\mathsf{CD}(K,N)$.
\begin{remark}
If a complete separable metric measure space $(\M,\mathsf{d},\m)$ satisfies the measure contraction property $\mathsf{MCP}(N-1,N)$ for some $N\ge1$ and if $\mathsf{diam}(M)\le\pi$ (which follows from the previous condition if $N\not=1$) then its $N$-Euclidean cone
$(\mathsf{Con}(\M),\mathsf{d_{Con}},\m_N)$  satisfies the measure contraction property $\mathsf{MCP}(0,N+1)$.
\end{remark}

As a second main result we deduce a generalized lower Ricci bound for the \textit{spherical cone} $\Sigma(\M)=\M\,\,{}_{\sin(r)}\!\!\times [0,\pi]$ over the compact Riemannian manifold $\M$. It can be defined as the quotient of the product space $\M\times[0,\pi]$ obtained by contracting all points in the fiber $\M\times\{0\}$ to the south pole $\mathcal{S}$ and all points in the fiber $\M\times\{\pi\}$ to the north pole $\mathcal{N}$. It is endowed with a metric $\mathsf{d}_\Sigma$ defined via
$$\cos\left(\mathsf{d}_\Sigma(p,q)\right)=\cos s\cos t+\sin s\sin t\cos\left(\mathsf{d}(x,y)\wedge \pi\right)$$
for $p=(x,s), q=(y,t)\in\Sigma(\M)$ and with a measure
$d\widehat{\m}_n(x,s):=d\mathsf{vol}(x)\otimes(\sin^ns ds)$.

\begin{theorem}
(i) The Ricci curvature of $\M$ is bounded from below by $n-1$ and $\mathsf{diam}(M)\le\pi$  if and only if the metric measure space $(\Sigma(\M),\mathsf{d_\Sigma},\widehat\m_n)$ satisfies the curvature-dimension condition $\mathsf{CD}(n,n+1)$.

(ii) A weighted Riemannian manifold satisfies the curvature-dimension condition $\mathsf{CD}(N-1,N)$ for a given real number $N>1$ if and only if the associated $N$-spherical cone satisfies the curvature-dimension condition $\mathsf{CD}(N,N+1)$.
\end{theorem}

Note that the analogous results holds true for generalized lower bounds for the sectional curvature.

\begin{remark}[see e.g. \cite{bi}, Theorem 4.7.1, 10.2.3] Let $(\M,\mathsf d)$ be a complete length metric  space with $\mathsf{diam}(M)\le\pi$.

(i) Then $(\M,\mathsf{d})$ has curvature bounded from below by 1 in the sense of Alexandrov if and only if the Euclidean cone $(\mathsf{Con}(\M),\mathsf{d_{Con}})$ has nonnegative curvature in the sense of Alexandrov.

(ii) Moreover, $(\M,\mathsf{d})$ has curvature bounded from below by 1 in the sense of Alexandrov if and only if the spherical cone $(\mathsf{Con}(\M),\mathsf{d_{Con}})$ has curvature bounded from below by 1 in the sense of Alexandrov.
\end{remark}
Note that the diameter bound is redundant if $M$ is not one-dimensional.

\medskip

Metric cones play an important role in the study of limits of Riemannian manifolds. Assume for instance that $(\M,\mathsf{d})$ is the Gromov-Hausdorff limit of a sequence of complete $n$-dimensional Riemannian manifolds whose Ricci curvature is uniformly bounded from below. Then in the non-collapsed case, every tangent cone $\mathsf{T}_x\M$ is a metric cone $\mathsf{Con}(\mathsf{S}_x\M)$ with $\mathsf{diam}(\mathsf{S}_x\M)\leq\pi$ \cite{cca,ccb}. The latter we would expect from the diameter estimate by Bonnet-Myers if $\mathsf{Ric}\geq n-2$ on $\mathsf{S}_x\M$ which in turn is consistent with the formal assertion \textquoteleft$\mathsf{Ric}\geq 0$ on $\mathsf{T}_x\M$\textquoteright.

\subsection{Basic Definitions and Notations}

Throughout this paper,  $(\mathsf{M},\mathsf{d})$ always will denote a complete separable metric space $(\mathsf{M},\mathsf{d})$ and $\mathsf{m}$ a locally finite measure  on $(\mathsf{M},\mathcal{B}(\mathsf{M}))$ with full support. That is,  for all $x\in \mathsf{M}$ and all sufficiently small $r>0$ the volume $\mathsf{m}(B_r(x))$ of balls centered at $x$ is positive and finite. To avoid pathologies, we assume that $M$ has more than one point.
Such a triple $(\mathsf{M},\mathsf{d},\m)$ will henceforth called \emph{metric measure space}.

The metric space $(\mathsf{M},\mathsf{d})$  is called a \textit{length space} iff $\mathsf{d}(x,y)=\inf\mathsf{Length}(\gamma)$ for all $x,y\in\mathsf{M}$, where the infimum runs over all curves $\gamma$ in $\mathsf{M}$ connecting $x$ and $y$. $(\mathsf{M},\mathsf{d})$ is called a \textit{geodesic space} if and only if every two points $x,y\in\mathsf{M}$ are connected by a curve $\gamma$ with $\mathsf{d}(x,y)=\mathsf{Length}(\gamma)$.
Distance minimizing curves of constant speed are called \textit{geodesics}.
The space of all  geodesics $\gamma:[0,1]\to M$ will be denoted by $\Gamma(M)$

$(\mathsf{M},\mathsf{d})$ is called \textit{non-branching} if for every tuple $(z,x_0,x_1,x_2)$ of points in $\mathsf{M}$ for which $z$ is a midpoint of $x_0$ and $x_1$ as well as of $x_0$ and $x_2$, it follows that $x_1=x_2$.

\medskip

$\mathcal{P}_2(\mathsf{M},\mathsf{d})$ denotes the \textit{$\mathsf{L}^2$-Wasserstein space} of probability measures $\mu$ on $(\mathsf{M},\mathcal{B}(\mathsf{M}))$ with finite second moments which means that $\int_\mathsf{M}\mathsf{d}^2(x_0,x)d\mu(x)<\infty$
for some (hence all) $x_0\in \mathsf{M}$. The \textit{$\mathsf{L}^2$-Wasserstein distance} $\mathsf{d}_{\mathsf{W}}(\mu_0,\mu_1)$ between two probability measures
$\mu_0,\mu_1\in\mathcal{P}_2(\mathsf{M},\mathsf{d})$ is defined as
$$\mathsf{d}_{\mathsf{W}}(\mu_0,\mu_1)=\inf\left\{\left(\int_{\mathsf{M}\times \mathsf{M}}\mathsf{d}^2(x,y)\,d\mathsf{q}(x,y)\right)^{1/2}:
\text{$\mathsf{q}$ coupling of $\mu_0$ and $\mu_1$}\right\}.$$
Here the infimum ranges over all \textit{couplings} of $\mu_0$ and $\mu_1$, i.e. over all probability measures on $\mathsf{M}\times \mathsf{M}$ with marginals $\mu_0$ and $\mu_1$. Equipped with this metric, $\mathcal{P}_2(\M,\mathsf{d})$ is a complete separable metric space.
The subspace of $\mathsf{m}$-absolutely continuous measures is denoted by $\mathcal{P}_2(\mathsf{M},\mathsf{d},\mathsf{m})$.

\medskip

\begin{definition}
\begin{itemize}
\item[(i)]
A subset $\Xi\subset\M\times\M$ is called \emph{$\mathsf{d}^2$-cyclically monotone} if and only if for any $k\in\N$ and for any family $(x_1,y_1),\dots,(x_k,y_k)$ of points in $\Xi$ the inequality
$$\sum^k_{i=1}\mathsf{d}^2(x_i,y_i)\leq\sum^k_{i=1}\mathsf{d}^2(x_i,y_{i+1})$$
holds with the convention $y_{k+1}=y_1$.
\item[(ii)] Given probability measures $\mu_0,\mu_1$ on $M$, a probability measure $\mathsf{q}$ on $M\times M$ is called
\emph{optimal coupling} of them iff $q$ has marginals $\mu_0$ and $\mu_1$ and
$$\mathsf{d^2_W}(\mu_0,\mu_1)=\int_{\M\times \M}\mathsf{d^2}(x,y) \, d\mathsf{q}(x,y).$$
\item[(iii)] A probability measure $\nu$ on $\Gamma(M)$ is called \emph{optimal path measure} (or dynamical optimal transference plan) iff the probability measure $(e_0,e_1)_*\nu$ on $M\times M$ is an optimal coupling of the probability measures
$(e_0)_*\nu$ and $(e_1)_*\nu$ on $M$.
\end{itemize}
\end{definition}
Here and in the sequel $e_t:\Gamma(M)\to M$ for $t\in[0,1]$ denotes the evaluation map $\gamma\mapsto \gamma_t$.
Moreover, for each measurable map $f:M\to M'$ and each measure $\mu$ on $M$ the  push forward (or image measure) of $\mu$ under $f$ will be denoted by $f_*\mu$.

\medskip

From  \cite[Lemma 2.11]{sa},\cite[Theorem 5.10]{vib} we quote:

\begin{lemma} \label{optimaltransport}
\begin{itemize}
\item[(i)] For each pair $\mu_0,\mu_1\in\mathcal{P}_2(\M,\mathsf{d})$ there exists an optimal coupling $\mathsf{q}$.
\item[(ii)]
The support of any  optimal coupling $\mathsf{q}$  is a $\mathsf{d}^2$-cyclically monotone set.
\item[(iii)] If $M$ is geodesic then
for each pair $\mu_0,\mu_1\in\mathcal{P}_2(\M,\mathsf{d})$ there exists an optimal path measure  with given initial and terminal distribution:
$(e_0)_*\nu=\mu_0$ and $(e_1)_*\nu=\mu_1$.
\item[(iv)] Given any optimal path measure $\nu$ as above, a geodesic $(\mu_t)_{t\in[0,1]}$ in
$\mathcal{P}_2(\M,\mathsf{d})$ connecting $\mu_0$ and $\mu_1$ is given by
$$\mu_t:=(e_t)_*\nu.$$
\item[(v)]If $(\M,\mathsf{d})$ is a non-branching space, then for each pair of geodesics $\gamma, \gamma'$  in the support of an optimal path measure we have:
 $$\gamma_{1/2}=\gamma'_{1/2} \quad\Longrightarrow\quad \gamma=\gamma'.$$
\end{itemize}
\end{lemma}

\subsection{The Curvature-Dimension Condition}

\begin{definition}
Given  $K\in\mathbb{R}$ and $N\in[1,\infty)$, the condition $\mathsf{CD}(K,N)$ states that for each pair
$\mu_0,\mu_1\in\mathcal{P}_2(\mathsf{M},\mathsf{d},\mathsf{m})$  there exist  an optimal
coupling $\mathsf{q}$ of $\mu_0=\rho_0\mathsf{m}$ and $\mu_1=\rho_1\mathsf{m}$ and a geodesic $\mu_t=\rho_t\, m$ in $\mathcal{P}_2(\mathsf{M},\mathsf{d},\mathsf{m})$ connecting them such that
\begin{equation} \label{CD}
\begin{split}
\int_\mathsf{M}\rho_t^{1-1/N'}d\m\ge\int_{\mathsf{M}\times
\mathsf{M}}\left[\tau^{(1-t)}_{K,N'}(\mathsf{d}(x_0,x_1))\rho^{-1/N'}_0(x_0)+
\tau^{(t)}_{K,N'}(\mathsf{d}(x_0,x_1))\rho^{-1/N'}_1(x_1)\right]d\mathsf{q}(x_0,x_1)
\end{split}
\end{equation}
for all $t\in (0,1)$ and all $N'\geq N$.
\end{definition}

In the case $K>0$, the \textit{volume distortion coefficients} $\tau^{(t)}_{K,N}(\cdot)$
for  $t\in (0,1)$  are defined by
$$\tau_{K,N}^{(t)}(\theta)=t^{1/N}\cdot
\left[ \frac {\sin\left(\sqrt{\frac K{N-1}}t\theta\right)} {\sin\left(\sqrt{\frac K{N-1}}\theta\right)} \right]^{1-1/N}$$
if $0\le\theta<\sqrt{\frac{N-1}K}\pi$ and by $\tau_{K,N}^{(t)}(\theta)=\infty$ if $\theta\ge\sqrt{\frac{N-1}K}\pi$.
In the case $K<0$ an analogous definition applies with $\sin\left(\sqrt{\frac K{N-1}}\ldots\right)$ replaced by
$\sinh\left(\sqrt{\frac {-K}{N-1}}\ldots\right)$. In the case $K=0$ simply
$$\tau_{0,N}^{(t)}(\theta)=t.$$
Therefore, the condition $\mathsf{CD}(0,N)$ just asserts that for each $N'\ge N$ the \emph{R\'enyi entropy}
$\mathsf{S}_{N'}(\nu_t|\mathsf{m}):=-\int_\mathsf{M}\rho_t^{1-1/N'}d\m$  is convex in $t\in[0,1]$.

 Replacing the volume distortion coefficients $\tau^{(t)} _{K,N}(\cdot)$ by slightly smaller coefficients $\sigma^{(t)} _{K,N}(\cdot)$ in the definition of $\mathsf{CD}(K,N)$ leads to the reduced curvature-dimension condition $\mathsf{CD}^*(K,N)$, a condition introduced and studied in \cite{bs},\cite{ds}.

The definitions of the condition $\mathsf{CD}(K,N)$ in \cite{sb} and \cite{lva} slightly differ. We follow the notation of \cite{sb}. For non-branching spaces, both concepts coincide.
In this case, it suffices to verify (\ref{CD}) for $N'=N$ since this already implies (\ref{CD}) for all $N'\ge N$.
Even more, the condition (\ref{CD}) can be formulated as a pointwise inequality.

\begin{lemma}[\cite{sb,lva,vib}]
 A nonbranching metric measure space $(\M,\mathsf{d},\m)$ satisfies the curvature dimension condition
$\mathsf{CD}(K,N)$ for  given  numbers $K$ and $N$ if and only if
for each pair
$\mu_0,\mu_1\in\mathcal{P}_2(\mathsf{M},\mathsf{d},\mathsf{m})$  there exist  an optimal path measure $\nu$
with initial and terminal distributions $(e_0)_\nu=\mu_0$, $(e_1)_\nu=\mu_1$ such that for $\nu$-a.e. $\gamma\in \Gamma(M)$ and all $t\in(0,1)$
\begin{equation}
\rho_t^{-1/N}(\gamma_t)
\ge\tau^{(1-t)}_{K,N}(\dot\gamma)\cdot\rho^{-1/N}_0(\gamma_0)+
\tau^{(t)}_{K,N}(\dot\gamma)\cdot\rho^{-1/N}_1(\gamma_1)
\end{equation}
where $\dot\gamma:=\mathsf{d}(\gamma_0,\gamma_1)$ and $\rho_t$ denotes the Radon-Nikodym density of $(e_t)_*\nu$ with respect to $\m$
\end{lemma}

\begin{lemma} \label{ohta} Assume that a metric measure space  $(\M,\mathsf{d},\m)$ satisfies the curvature dimension condition
$\mathsf{CD}(N-1,N)$ for some number $N>1$.
\begin{itemize}
\item[(i)] Then the diameter of $M$ is bounded by $\pi$.

\item[(ii)] Moreover, for every $x\in\M$ the set $\mathsf{M}_x:=\{x'\in\M:\mathsf{d}(x,x')=\pi\}$ of antipodes of $x$ consists of at most one point.
\end{itemize}
\end{lemma}
Assertion (i), the 'generalized Bonnet-Myers theorem' was proven in \cite{sb}. Assertion (ii) is due to Shin-ichi Ohta \cite[Theorem 4.5]{oa}.

\bigskip

Now let us have closer look on the curvature-dimension condition in the case of weighted Riemannian spaces.
Given a complete $n$-dimensional manifold $M$ equipped with its Riemannian distance $\mathsf d$ and with a weighted measure $d\m(x)=e^{-V(x)}d\mathsf{vol}_\M(x)$ for some function $V:\M\to\R$.  Then for each real number $N>1$
the $N$-Ricci tensor is defined as
$$\mathsf{Ric}^{N,V}_x(v,v):=
\mathsf{Ric}_x(v,v)
+\left[ \mathsf{Hess}\, V - \frac1{N-n}\nabla V \otimes \nabla V\right]_{x} (v,v).$$

\begin{lemma}[\cite{lva,sb}]
The weighted Riemannian space $(\M,\mathsf{d},\m)$ satisfies the condition $\mathsf{CD}(K,N)$ if and only if
$\mathsf{Ric}^{N,V}\ge K$ on $M$ in the sense that
$$\mathsf{Ric}^{N,V}_x(v,v)\ge K\cdot\|v\|_{T_x}^2$$
for all $x\in M$ and all $v\in T_xM$.
\end{lemma}

\section{Euclidean Cones over Metric Measure Spaces} \label{cones}

\begin{definition}[$N$-Euclidean cone]
For a metric measure space $(\mathsf{M},\mathsf{d},\mathsf{m})$ and any  $N\in [1,\infty)$, the \textit{$N$-Euclidean cone} $(\mathsf{Con}(\M),\mathsf{d_{Con}},\m_N)$ is a metric measure space defined as follows:
\begin{itemize}
\item[$\diamond$] $\mathsf{Con}(\M):=\mathsf{M}\times[0,\infty)\ / \ \mathsf{M}\times\{0\}$
\item[$\diamond$] For $(x,s),(x',t)\in \mathsf{M}\times[0,\infty)$
$$\mathsf{d_{Con}}((x,s),(x',t)):=\sqrt{s^2+t^2-2st\cos\left(\mathsf{d}(x,x')\wedge\pi\right)}$$
\item[$\diamond$] $d\m_N(x,s):=d\mathsf{m}(x)\otimes s^Nds$.
\end{itemize}
The point $\mathsf{O}:=\M\times\{0\}\in\mathsf{Con}(\M)$ is called origin of the cone.
\end{definition}

The most prominent example in this setting is the unit sphere $\s^n\subset\R^{n+1}$, endowed with its intrinsic Riemannian distance and with the Riemannian volume measure  on it. In other words, $\mathsf d(x,y)$ is the Euclidean angle between the rays from the origin $0\in\R^{n+1}$ to the points $x$ and $y$ on the unit sphere of $\R^{n+1}$. Each $\xi\in\R^{n+1}\setminus\{0\}$ can be uniquely written as $\xi=(x,r)$ with $r\in(0,\infty)$ and $x\in \s^n$, namely, $r=|\xi|$ and $x=\frac\xi{|\xi|}$.

The definition of the metric $\mathsf{d_{Con}}$ and the measure $\m_n$ ensures that the $n$-Euclidean cone over $\s^n$ is the Euclidean space $\R^{n+1}$ equipped with the Euclidean metric and the Lebesgue measure expressed in spherical coordinates.

\begin{conjecture}
A metric measure space $(\M,\mathsf{d},\m)$ satisfies the curvature-dimension condition $\mathsf{CD}(N-1,N)$ for some real $N\ge1$ and  $\mathsf{diam}(M)\le\pi$ (which follows from the previous condition if $N\not=1$) if and only if the $N$-Euclidean cone $(\mathsf{Con}(\M),\mathsf{d_{Con}},\m_N)$ satisfies the curvature-dimension condition $\mathsf{CD}(0,N+1)$.
\end{conjecture}

Conjecture 2.2 is true for every weighted Riemannian space.
 The proof  is based on two ingredients:
\begin{itemize}
\item[(a)] Optimal transports on the cone never transport mass through the origin, -- provided the base space $M$ satisfies an appropriate $\mathsf{CD}$ condition.
\item[(b)] Optimal transports on the punctured cone $\mathsf{Con}_0(\M)$ satisfy the $\mathsf{CD}$ condition implied by the Ricci bound for the non-complete, weighted Riemannian manifold $\mathsf{Con}_0(\M)$.
The latter in turn is equivalent to a Ricci bound for the complete weighted Riemannian manifold $\M$.
\end{itemize}
Property (a) will be proven  as a result of independent interest for general metric measure spaces.

\medskip

\begin{theorem}\label{origin}
Assume that the metric measure space $(\mathsf{M},\mathsf{d},\mathsf{m})$ satisfies the curvature-dimension condition $\mathsf{CD}(N-1,N)$ for some  $N\ge1$ and that $\mathsf{diam}(M)\le\pi$ (which follows from the previous condition if $N\not=1$). Let $\nu$ be any optimal path measure on the Euclidean cone $(\mathsf{Con}(\M),\mathsf{d_{Con}})$.
\begin{itemize}
\item[(i)] For every $t\in(0,1)$ there exists at most one geodesic $\gamma\in\mathsf{supp}[\nu]$ with $\gamma_t =\mathsf O$.
\item[(ii)]
For every $r>0$ there exists at most one $x\in M$ such that $\gamma_0=(x,r)$ is the initial point of some geodesic $\gamma\in \mathsf{supp}[\nu]\cap \Gamma_{\mathsf O}$ where $$\Gamma_{\mathsf O}:=\{\gamma\in\Gamma(\mathsf{Con}(\M)): \  \gamma_t=\mathsf O\ \mbox{\rm for some }t\in(0,1)\}.$$
\item[(iii)]  If $(e_0)_*\nu\ll m_N$ then $\nu$ gives no mass to geodesics through $\mathsf O$:
$$\nu\left(\Gamma_{\mathsf O}\right)=0.$$
\end{itemize}
\end{theorem}
\vspace{-0.8cm}
\begin{center}
\hspace{-2cm}
\includegraphics[height=4.5cm]{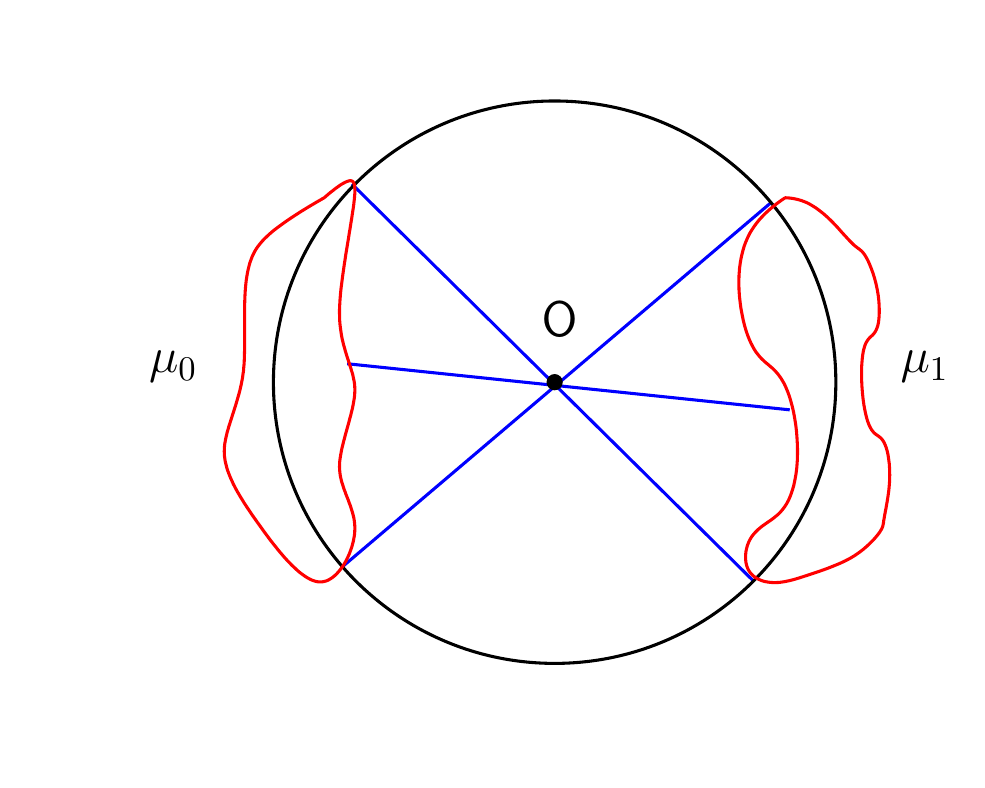}
\end{center}
\vspace{-0.5cm}
\begin{proof}
(i) Fix $t\in(0,1)$ and assume that two geodesics $\gamma,\gamma'\in\mathsf{supp}[\nu]$ have the origin as common $t$-intermediate point, i.e. $\gamma_t=\gamma'_t=\mathsf O$. Then $\gamma_0=(x_0,tr)$, $\gamma_1=(x_1,(1-t)r)$ for some $x_0,x_1\in M$ and with $r=\dot\gamma=\mathsf d_{Con}(\gamma_0,\gamma_1)$.
Similarly,  $\gamma'_0=(x'_0,tr')$, $\gamma'_1=(x'_1,(1-t)r')$ for some $x'_0,x'_1\in M$ and with $r'=\dot\gamma'$. If $r>0$ then $x_0$ and $x_1$  are antipodes of each other (i.e. $\mathsf d(x_0,x_1)=\pi$). Similarly, for $x'_0$ and $x_1'$.  (See e.g. Lemma \ref{antispher} for a detailed proof in the more sophisticated case of spherical cones.)

Cyclic monotonicity implies
$$0\le \mathsf d^2_{Con}(\gamma_0,\gamma'_1)+\mathsf d^2_{Con}(\gamma'_0,\gamma_1)-
\mathsf d^2_{Con}(\gamma_0,\gamma_1)-\mathsf d^2_{Con}(\gamma'_0,\gamma'_1).$$
On the other hand, a simple application of the triangle inequality yields
\begin{eqnarray*}
\mathsf d^2_{Con}(\gamma_0,\gamma'_1)+\mathsf d^2_{Con}(\gamma'_0,\gamma_1)&-&
\mathsf d^2_{Con}(\gamma_0,\gamma_1)-\mathsf d^2_{Con}(\gamma'_0,\gamma'_1)\\
&\le&\left[tr+(1-t)r'\right]^2+\left[tr'+(1-t)r\right]^2-r^2-r'^2\\
&=&-2t(1-t)(r-r')^2.
\end{eqnarray*}
Hence, $r=r'$.

With this at hand, a more precise calculation yields
\begin{eqnarray*}
0&\le&\mathsf d^2_{Con}(\gamma_0,\gamma'_1)+\mathsf d^2_{Con}(\gamma'_0,\gamma_1)-
\mathsf d^2_{Con}(\gamma_0,\gamma_1)-\mathsf d^2_{Con}(\gamma'_0,\gamma'_1)\\
&=& 2r^2\left[t^2+(1-t)^2-t(1-t)\cos\mathsf d(x_0,x_1')-t(1-t)\cos\mathsf d(x'_0,x_1)\right]-2r^2\\
&=& -2r^2t(1-t)\left[2+\cos\mathsf d(x_0,x_1')+\cos\mathsf d(x'_0,x_1)\right].
\end{eqnarray*}
Thus $\mathsf d(x_0,x_1')=\mathsf d(x'_0,x_1)=\pi$. That is, $x_0$ and $x'_1$ are antipodes (as well as $x'_0$ and $x_1$). Since antipodes in $M$ are unique (Lemma \ref{ohta}(ii)) we conclude that $x_0=x'_0$ and $x_1=x'_1$.
Thus $\gamma_0=\gamma'_0$ and $\gamma_1=\gamma'_1$.

\medskip

In most cases of interest, geodesics are uniquely determined by their initial and terminal points. In these case, we are done. The general case, requires an additional argument.
An optimal path measure $\nu$ not only induces an optimal coupling $(e_0,e_1)_*\nu$ between its initial and terminal distribution $(e_0)_*\nu$ and $(e_1)_*\nu$, resp. More generally,  the measure
$(e_\sigma,e_\tau)_*\nu$ will be an optimal coupling of $(e_\sigma)_*\nu$ and $(e_\tau)_*\nu$ for each $0\le\sigma\le\tau\le1$.
For each $\sigma\in(0,t)$ one can choose $\tau\in(t,1)$ (and vice versa) such that  $(e_t)_*\nu$ is a
 $t$-intermediate point of $(e_\sigma)_*\nu$ and $(e_\tau)_*\nu$. Hence, the previous argument will imply that
$\gamma_\sigma=\gamma'_\sigma$ and $\gamma_\tau=\gamma'_\tau$.
This finishes the proof.

\medskip

(ii)
Assume $\gamma,\gamma'\in \mathsf{supp}[\nu]\cap \Gamma_{\mathsf O}$ with $\gamma_0=(x_0,r)$ and $\gamma'_0=(x'_0,r)$. The fact that $\gamma$ passes through the origin implies that $\gamma_1=(x_1,r_1)$ with $x_1\in M$ being an antipode of $x_0$, i.e. $\mathsf{d}(x_0,x_1)=\pi$.
Similarly,  $\gamma'_1=(x'_1,r'_1)$ with  $\mathsf{d}(x'_0,x'_1)=\pi$. The radii $r_1, r'_1$ are arbitrary positive numbers. Cyclic monotonicity implies
 \begin{eqnarray*}
0&\le&\mathsf d^2_{Con}(\gamma_0,\gamma'_1)+\mathsf d^2_{Con}(\gamma'_0,\gamma_1)-
\mathsf d^2_{Con}(\gamma_0,\gamma_1)-\mathsf d^2_{Con}(\gamma'_0,\gamma'_1)\\
&=& r^2+r_1'^2 -2rr_1'\cos\mathsf d(x_0,x_1')+r^2+r_1^2-2rr_1\cos\mathsf d(x'_0,x_1)-(r+r_1)^2-(r+r'_1)^2\\
&=& -2rr_1'\left[1+\cos\mathsf d(x_0,x_1')\right]-2rr_1\left[1+\cos\mathsf d(x'_0,x_1)\right].
\end{eqnarray*}
Hence, $\mathsf d(x_0,x_1')=\pi$. That is, $x_0$ and $x'_1$ are antipodes (as well as $x'_0$ and $x'_1$, which has been observed before). Uniqueness of antipodes in $M$ implies $x_0=x'_0$.

\medskip

(iii) Let us assume that $\nu(\Gamma_{\mathsf O})>0$. Then without restriction we even may assume that $\nu$ is supported by $\Gamma_{\mathsf O}$. (Otherwise, replace $\nu$ by its restriction onto the set $\Gamma_{\mathsf O}$.) Since $\m_N({\mathsf O})=0$ we may also assume that $\gamma_0\not={\mathsf O}$ and $\gamma_1\not={\mathsf O}$ for $\nu$-a.e. $\gamma$.

The previous part (ii) asserts that for each $r>0$ there exists at most one point $x_0=f(r)\in M$ such that $(f(r),r)$ is the initial point $\gamma_0$ of some geodesic $\gamma\in\mathsf{supp}[\nu]\cap \Gamma_{\mathsf O}$.
Thus the measure $\mu_0:=(e_0)_*\nu$ is concentrated on the set
 $C_f:=\{(f(r),r)\in \mathsf{Con}(\M):\, r>0\}$.

The curvature-dimension condition for the base space $M$ implies that $\m$ has no atoms. Hence,
$$\m_N(C_f)=0$$
and therefore $\mu_0\not\ll \m_N$.
\end{proof}

According to the previous result, we know that  -- under the given curvature-dimension assumptions -- optimal path measures on an Euclidean cone never will transport mass through the origin. It therefore suffices to study optimal transports on the
 \emph{punctured cone}
$$\mathsf{C}_0:=\mathsf{Con}(\M)\setminus\{\mathsf{O}\}.$$
To analyze such transports,
 we restrict ourselves to base spaces $M$ which are (weighted) Riemannian manifolds. Our results
crucially will rely on the fact that in this case the punctured cone
$\mathsf{C}_0$
is a non-complete(!) Riemannian manifold and that the Ricci curvature of it can be calculated explicitly.
More precisely, the punctured $n$-Euclidean cone is a Riemannian manifold whereas the punctured $N$-Euclidean cone is a weighted Riemannian manifold.

\begin{lemma}
(i) The punctured Euclidean cone $\mathsf{C}_0$ is an $(n+1)$-dimensional Riemannian manifold. For $(x,r)\in\mathsf{C}_0$ with $x\in\M$ and $r>0$ the tangent space $T_{(x,r)}\mathsf{C}_0$ can be parametrized as $T_x\M\oplus\R$ with $\|(v,t)\|^2_{T_{(x,r)}}=r^2\,\|v\|^2_{T_x}+t^2$.
Moreover, for $(v,t)\in{T}_{(x,r)}\mathsf{C}_0$ with $v\in{T}_x\M$ and $t\in\R$ we have the identity
$$\mathsf{Ric}_{(x,r)}((v,t),\, (v,t))=
\mathsf{Ric}_x(v,v)-(n-1)\| v\|^2_{T_x}.$$
In particular,  $\mathsf{Ric}\geq 0$  on $\mathsf{C}_0$ if and only if $\mathsf{Ric}\geq n-1$ on $\M$.

(ii) The punctured $N$-Euclidean cone $\mathsf{C}_0$ is a weighted $(n+1)$-dimensional Riemannian manifold with measure
$d\m_N(x,r)=r^N\,dr\, d\mathsf{vol}_\M(x)=e^{-W(r)}\,d \m_n(x,r)$ where $W(r)=-(N-n)\log r$ and $d\m_n(x,r)=d\mathsf{vol}_{\mathsf{C}_0}(x,r)=r^n\,dr\, d\mathsf{vol}_\M(x)$ denotes the Riemannian volume measure on $\mathsf{C}_0$.
For each $N\ge n$, the $(N+1)$-Ricci tensor satisfies
$$\mathsf{Ric}^{N+1,W}_{(x,r)}((v,t),\, (v,t))=
\mathsf{Ric}_x(v,v)-(N-1)\| v\|^2_{T_x}.$$
In particular,  $\mathsf{Ric}^{N+1,W}\geq 0$ on $\mathsf{C}_0$ if and only if $\mathsf{Ric}\geq N-1$ on $\M$.

(iii) More generally, let the $n$-dimensional Riemannian manifold $\M$ be equipped with the weighted measure
$d\m(x)=e^{-V(x)}\,d\mathsf{vol}_\M(x)$ for some $V:\M\to\R$ and let the punctured cone $\mathsf{C}_0$ be equipped with the measure
$$d\m_N(x,r)=r^Ndr\, d\m(x)=e^{-V(x)-W(r)}\, d\mathsf{vol}_{\mathsf{C}_0}(x,r)$$
with (as before) $W(r)=-(N-n)\log r$ and $d\mathsf{vol}_{\mathsf{C}_0}(x,r)=r^n\,dr\, d\mathsf{vol}_\M(x)$. Then
\begin{equation}\mathsf{Ric}^{N+1,V+W}_{(x,r)}((v,t),\, (v,t))=
\mathsf{Ric}^{N,V}_x(v,v)-(N-1)\| v\|^2_{T_x}.
\end{equation}
In particular,  $\mathsf{Ric}^{N+1,V+W}\geq 0$ on $\mathsf{C}_0$ if and only if $\mathsf{Ric}^{N,V}\geq N-1$ on $\M$.
\end{lemma}

\begin{proof}
Assertion (i) is a classical result due to  Cheeger and Taylor \cite{cha, chb}. Assertion (ii) is the particular case of (iii) with $V=0$.

\medskip

(iii):  For arbitrary $V(x,r)=V(x)$ and $W(x,r)=W(r)$ depending only on the radial coordinate $r\in\R$ or on the basic coordinate $x\in M$, resp., we have
$$\nabla W_{(x,r)}(v,t)=h'(0),\qquad
\left[  \mathsf{Hess}\, W\right]_{(x,r)} ((v,t),\, (v,t))=h''(0)
$$
for all $(x,r)\in \mathsf{C}_0$ and all $(v,t)\in T_{(x,r)}\mathsf{C}_0$
where
$h(s):=W\left(\sqrt{(r+st)^2+s^2r^2\| v\|^2_{T_x}}\right)$.
Moreover,
$$\left[ \nabla V \otimes \nabla W\right]_{(x,r)} ((v,t),\, (v,t))\quad=\quad{\nabla V}_x(v)\cdot W'(r)\cdot t$$ for all $(x,r)\in \mathsf{C}_0$ and all $(v,t)\in T_{(x,r)}\mathsf{C}_0$ as well as
$$
\left[  \mathsf{Hess}\, V\right]_{(x,r)} ((v,t),\, (v,t))\quad=\quad f''(0)\quad=\quad\left[  \mathsf{Hess}\, V\right]_{x} (v,v)-2{\nabla V}_x(v)\cdot \frac t r$$
where the expressions on the LHS always have to be interpreted as quantities on the $(n+1)$-dimensional manifold $ \mathsf{C}_0$ whereas the expressions on the RHS are the original data on the basic $n$-dimensional manifold $M$ and where
$$f(s)=V\left(\exp_x\left(\frac v{\|v\|_{T_x}}\cdot\arctan\frac{rs \|v\|_{T_x}}{r+st}\right)\right).$$
For the particular choice of $W(x,r)=-(N-n)\log r$ explicit calculations yield
$$\left[  \mathsf{Hess}\, W - \frac1{N-n}\nabla W \otimes \nabla W\right]_{(x,r)} ((v,t),\, (v,t))=-(N-n)\, \| v\|^2_{T_x}.$$
Hence, together with the identity from (i)
\begin{eqnarray*}
\lefteqn{\mathsf{Ric}^{N+1, V+W}_{(x,r)}((v,t),\, (v,t))}\\
&=&
\mathsf{Ric}_{(x,r)}((v,t),\, (v,t))
+\left[  \mathsf{Hess}\, (V+W) - \frac1{N-n}\nabla (V+W) \otimes \nabla (V+W)\right]_{(x,r)} ((v,t),\, (v,t))\\
&=&
\mathsf{Ric}_x(v,v)-(n-1)\| v\|^2_{T_x}
+\left[ \mathsf{Hess}\, V - \frac1{N-n}\nabla V \otimes \nabla V\right]_{x} (v,v)
-(N-n)\, \| v\|^2_{T_x}\\
&=&
\mathsf{Ric}^{N,V}_x(v,v)-(N-1)\| v\|^2_{T_x}.
\end{eqnarray*}

\end{proof}

\begin{theorem}
Given a complete $n$-dimensional manifold $M$ equipped with its Riemannian distance $\mathsf d$ and with a weighted measure $d\m(x)=e^{-V(x)}d\mathsf{vol}_\M(x)$ for some function $V:\M\to\R$.  Then for each real number $N>1$ the following statements are equivalent:
\begin{itemize}
\item[(i)]
The weighted Riemannian space $(\M,\mathsf{d},\m)$ satisfies the  condition $\mathsf{CD}(N-1,N)$.
\item[(ii)] The $N$-Euclidean cone $(\mathsf{Con}(\M),\mathsf{d_{Con}},\m_N)$  satisfies the condition $\mathsf{CD}(0,N+1)$.
\end{itemize}
\end{theorem}

\begin{proof}
$(ii) \Rightarrow (i)$:\ The curvature-dimension condition  $\mathsf{CD}(0,N+1)$ for the $N$-Euclidean cone $(\mathsf{Con}(\M),\mathsf{d_{Con}},\m_N)$  implies that this condition holds \emph{locally} on the punctured cone. For this (non-complete) weighted Riemannian manifold, however, the local curvature-dimension condition $\mathsf{CD}_{loc}(0,N+1)$ is equivalent to nonnegativity of the $(N+1)$-Ricci tensor $\mathsf{Ric}^{N+1,V+W}$ on $\mathsf{C}_0$, see Lemma 1.11. Due to the previous Lemma 2.4(iii), this implies $\mathsf{Ric}_M^{N,V}\ge N-1$. For the (complete) weighted Riemannian space  $(\M,\mathsf{d},\m)$, the latter in turn is equivalent to
$\mathsf{CD}(N-1,N)$.

$(i) \Rightarrow (ii)$:\  Let probability measures $\mu_0$ and $\mu_1$ on $\mathsf{Con}(\M)$ be given, absolutely continuous with respect to $\m_N$. According to Theorem \ref{origin}, any optimal path measure $\nu$ with marginal distributions $(e_0)_*\nu=\mu_0$ and $(e_1)_*\nu=\mu_1$ will give no mass to geodesics through the origin. In other words, $\nu$-almost every geodesic will  stay within the punctured cone $\mathsf{C}_0$.

According to  Lemma 2.4(iii), assertion (i) implies that the $(N+1)$-Ricci tensor $\mathsf{Ric}^{N+1,V+W}$ on the weighted Riemannian space $\mathsf{C}_0$ is nonnegative.
Hence, classical arguments based on Jacobi field calculus -- exactly the same as used  to deduce Lemma 1.11 -- will imply that (1.2) holds true with $K=0$ for $\nu$-a.e. geodesic $\gamma$ which remains within $\mathsf{C}_0$. That is, $\mathsf{CD}(0,N+1)$ holds true on $\mathsf{Con}(\M)$.
\end{proof}

\begin{corollary}
Given a complete $n$-dimensional manifold $M$ (equipped with its Riemannian distance $\mathsf d$ and its Riemannian volume $d\m=d\mathsf{vol}_M$) and a real number $N\ge1$. Then the following statements are equivalent:
\begin{itemize}
\item[(i)] $\mathsf{Ric}\ge N-1$ on $M$, $\mathsf{dim}_M\le N$ and $\mathsf{diam}(M)\le\pi$ (the latter follows from the Ricci and dimension bounds if $N\not=1$);
\item[(ii)]
The space $(\M,\mathsf{d},\m)$ satisfies the curvature-dimension condition $\mathsf{CD}(N-1,N)$ and $\mathsf{diam}(M)\le\pi$ (which follows from the $\mathsf{CD}$ condition if $N\not=1$);
\item[(iii)] The $N$-Euclidean cone $(\mathsf{Con}(\M),\mathsf{d_{Con}},\m_N)$  satisfies $\mathsf{CD}(0,N+1)$.
\end{itemize}
\end{corollary}

\begin{proof}
The equivalence $(i) \Leftrightarrow (ii)$ is well-known. Moreover, it is well-known that for each $N>1$ the condition  $\mathsf{CD}(N-1,N)$ implies $\mathsf{diam}(M)\le\pi$. See Lemma 1.10 and 1.11.

In the case $N\not=1$, the equivalence $(ii) \Leftrightarrow (iii)$  for Riemannian spaces  follows from the more general assertion of Theorem 2.5  for \emph{weighted} Riemannian spaces. Indeed, the arguments there also apply to the case $N=1$. It only remains to prove that (iii) in the case $N=1$ implies $\mathsf{diam}(M)\le\pi$.

\medskip

Assume the contrary: i.e. $M$ is a circle or an interval of $\mathsf{diam}(M)>\pi$. Then there exist non-empty intervals $I,J\subset M$ of length $R>0$ such that $d(x,y)>\pi$ for all $x\in I, y\in J$. Thus for all $x\in I,y\in J$ and $r\in(0,\infty)$ the origin $\mathsf O$ will be the unique midpoint of $(x,r)$ and $(y,r)$ in $\mathsf{Con}(\M)$. Moreover, for each pair $(x,s)\in I\times [1-\epsilon,1]$ and $(y,t)\in J\times [1-\epsilon,1]$ in $\mathsf{Con}(\M)$ the midpoint will lie in $B_\epsilon:=\left((I\cup J)\times(0,\epsilon]\right)\cup\{\mathsf O\}$.

Let $\mu_0$ and $\mu_1$ be the 'uniform distributions' on $I_\epsilon:=I\times [1-\epsilon,1]$ and $J_\epsilon:=J\times [1-\epsilon,1]$, resp., i.e. $d\mu_0=C_\epsilon\,1_{I_\epsilon}\,d\m_N$, $d\mu_1=C_\epsilon\,1_{J_\epsilon}\,d\m_N$ with suitable $C_\epsilon\ge\frac1{R\epsilon}$.
Then the Renyi entropy of them satisfies
$$-S_{N+1}(\mu_0|\m_N)=-S_{N+1}(\mu_1|\m_N)={C_\epsilon}^{-\frac1{N+1}}\le (R\epsilon)^{\frac1{N+1}}=c\, \epsilon^{\frac1{N+1}}.$$
On the other hand, the midpoint $\mu_{1/2}$ of $\mu_0$ and $\mu_1$ is supported on $B_\epsilon$. Hence, its Renyi entropy is bounded from below by the Renyi entropy of the uniform distribution on $B_\epsilon$:
$$S_{N+1}(\mu_{1/2}|\m_N)\ge S_{N+1}(C'_\epsilon\,1_{B_\epsilon}\, \m_N|\m_N)=-{C'_\epsilon}^{-\frac1{N+1}}=-c'\,\epsilon.$$
Note that ${C'_\epsilon}^{-1}=m_N(B_\epsilon)=2R\cdot\int_0^\epsilon r^N\,dr=\frac{2R}{N+1}\epsilon^{N+1}$.
Thus, choosing  $\epsilon$ sufficiently small we obtain
$$S_{N+1}(\mu_{1/2}|\m_N) \gg \frac12(S_{N+1}(\mu_0|\m_N)+S_{N+1}(\mu_1|\m_N))$$
which contradicts the $\mathsf{CD}(0,N+1)$ condition.
\end{proof}

\begin{example}
Let $\M=\left(\frac1{\sqrt 3} \s^2\right)\times\left(\frac1{\sqrt 3} \s^2\right)$.
\begin{itemize}
\item[(i)] Then the Euclidean cone over $\M$ -- more precisely, the metric measure space $(\mathsf{Con}(\M),\mathsf{d_{Con}},\m_4)$ -- satisfies the curvature-dimension condition
$\mathsf{CD}(0,5)$.
\item[(ii)] On the other hand, the Euclidean cone over $\M$ -- more precisely, the metric space $(\mathsf{Con}(\M),\mathsf{d_{Con}})$ -- is not an Alexandrov space: the sectional curvature on the punctured cone $\mathsf{C}_0$ is unbounded from below (and above) in any punctured neighborhood of the origin $\mathsf 0$.
\end{itemize}
\end{example}

\begin{proof}
Given $x,y\in\frac1{\sqrt 3} \s^2$, let $u_1,u_2$ be an orthonormal basis of $T_x(\frac1{\sqrt 3} \s^2)$ and
$v_1,v_2$ be an orthonormal basis of $T_y(\frac1{\sqrt 3} \s^2)$. Then an orthonormal basis of $T_{(x,y)}\M=T_x(\frac1{\sqrt 3} \s^2)\oplus T_y(\frac1{\sqrt 3} \s^2)$
is given by $\{\tilde u_1,\tilde u_2,\tilde v_1,\tilde v_2\}$ with $\tilde u_i=(u_i,0)$ and $\tilde v_i=(0,v_i)$.
In this basis
$$\mathsf{Sec}_{(x,y)}(\tilde u_1,\tilde u_2)=3, \quad
\mathsf{Sec}_{(x,y)}(\tilde u_1,\tilde v_1)=0,\quad \mathsf{Sec}_{(x,y)}(\tilde u_1,\tilde v_2)=0$$
and analogously for any other basis vector in the place of $\tilde u_1$.
Hence, in particular,
$$\mathsf{Ric}_{(x,y)}(\xi,\xi)=3$$
for each $\xi\in\{\tilde u_1,\tilde u_2,\tilde v_1,\tilde v_2\}$ and thus for each $\xi\in T_{(x,y)}\M$.

\smallskip

(i) Thus according to  Theorem 2.5 the Euclidean cone satisfies the $\mathsf{CD}(0,5)$ condition.

\smallskip

(ii) Given $r>0$ an orthonormal basis of $T_{(x,y,r)}\mathsf{C}_0=T_x(\frac1{\sqrt 3} \s^2)\oplus T_y(\frac1{\sqrt 3} \s^2)\oplus\R$ is given by $\{\hat u_1,\hat u_2,\hat v_1,\hat v_2,\hat w\}$ with $\hat u_i=\frac1r(u_i,0,0)$, $\hat v_i=\frac1r(0,v_i,0)$ and $\hat w=(0,0,1)$.
In this basis
\begin{equation}\mathsf{Sec}_{(x,y,r)}(\hat u_1,\hat u_2)=\frac2{r^2}, \quad\mathsf{Sec}_{(x,y,r)}(\hat u_1,\hat v_1)=-\frac1{r^2},\quad
\mathsf{Sec}_{(x,y,r)}(\hat u_1,\hat v_2)=-\frac1{r^2},\quad \mathsf{Sec}_{(x,y,r)}(\hat u_1,\hat w)=0\end{equation}
and analogously for $\tilde u_2, \tilde v_1$ or $\tilde v_2$  in the place of $\tilde u_1$.
Of course, this in particular implies
$\mathsf{Ric}_{(x,y,r)}(\xi,\xi)=0$
for each $\xi\in T_{(x,y,r)}\mathsf{C}_0$, -- in accordance with Lemma 2.4.
\end{proof}

\section{Spherical Cones over Metric Measure Spaces} \label{spherical}

There are further objects with famous Euclidean ancestors -- among them is the \textit{spherical cone} or \textit{suspension} over a topological space $\M$. We begin with a familiar example: In order to construct the Euclidean sphere $\s^{n+1}$ out of its equator $\s^n$ we add two poles $\mathcal{S}$ and $\mathcal{N}$ and connect them via semicircles, the \textit{meridians}, through every point in $\s^n$.

In the general case of abstract spaces $\M$, we consider the product $\M\times [0,\pi]$ and contract each of the fibers $\mathcal{S}:=\M\times\{0\}$ and $\mathcal{S}:=\M\times\{\pi\}$ to a point, the \textit{south} and the \textit{north pole}, respectively. The resulting space is denoted by $\Sigma(\M)$ and is called the spherical cone over $\M$.

\begin{definition}[$N$-spherical cone]
The \textit{$N$-spherical cone} $(\Sigma(\M),\mathsf{d}_\Sigma,\widehat\m_N)$ over a metric measure space $(\mathsf{M},\mathsf{d},\mathsf{m})$  is the metric measure space defined as follows:
\begin{itemize}
\item[$\diamond$] $\Sigma(\M):=\mathsf{M}\times[0,\pi]\ \Big/ \ \mathsf{M}\times\{0\},\M\times\{\pi\}$
\item[$\diamond$] For $(x,s),(x',t)\in \mathsf{M}\times[0,\pi]$
$$\cos\left(\mathsf{d}_\Sigma((x,s),(x',t))\right):=\cos s\cos t+\sin s\sin t\cos\left(\mathsf{d}(x,x')\wedge\pi\right)$$
\item[$\diamond$] $d\widehat\m_N(x,s):=d\m(x)\otimes(\sin^Ns ds)$.
\end{itemize}
\end{definition}

For a nice introduction and detailed information about Euclidean and spherical cones over metric spaces we refer to \cite{bi}.

\begin{lemma} \label{antispher}
Assume that $\mathsf{diam}(M)\le\pi$. Let  $\gamma:[0,1]\to\Sigma(\M)$ be a non-constant geodesic with endpoints $\gamma_0=(x_0,r_0)$ and $\gamma_1=(x_1,r_1)$ in $\Sigma(\M)$. If $\gamma_t=\mathcal{S}$ for some $t\in(0,1)$, then $x_0$ and $x_1$ are antipodes in $\M$.
\end{lemma}
\begin{proof} Due to the definition of $\mathsf{d}_\Sigma$, it holds that
$r_0=\mathsf{d}_\Sigma(\gamma_0,\gamma_t)=t\mathsf{d}_\Sigma(\gamma_0,\gamma_1)$
as well as
$r_1=\mathsf{d}_\Sigma(\gamma_t,\gamma_1)=(1-t)\mathsf{d}_\Sigma(\gamma_0,\gamma_1)$
and consequently, $r_1=\tfrac{1-t}{t}r_0$. Inserting this equality in the expression for $\cos\left(\tfrac{r_0}{t}\right)$ we obtain
\begin{align*}
\cos\left(\tfrac{r_0}{t}\right)&=\cos\left(\mathsf{d}_\Sigma(\gamma_0,\gamma_1)\right)=\cos r_0\cos\left(\tfrac{1-t}{t}r_0\right)+\sin r_0\sin\left(\tfrac{1-t}{t}r_0\right)\cos\left(\mathsf{d}(x_0,x_1)\right).
\end{align*}
Since $\mathsf{diam}(M)\le\pi$ by assumption, this leads to
\begin{align*}
\cos(\mathsf{d}(x_0,x_1))&=\frac{\cos\left(\tfrac{r_0}{t}\right)-\cos r_0\cos\left(\tfrac{1-t}{t}r_0\right)}{\sin r_0\sin\left(\tfrac{1-t}{t}r_0\right)}\\
&=\frac{\cos\left(\tfrac{r_0}{t}\right)-\tfrac{1}{2}\left[\cos\left(\tfrac{2t-1}{t}r_0\right)+
\cos\left(\tfrac{r_0}{t}\right)\right]}
{\tfrac{1}{2}\left[\cos\left(\tfrac{2t-1}{t}r_0\right)-\cos\left(\tfrac{r_0}{t}\right)\right]}\\
&=\frac{\tfrac{1}{2}\left[\cos\left(\tfrac{r_0}{t}\right)-\cos\left(\tfrac{2t-1}{t}r_0\right)\right]}
{\tfrac{1}{2}\left[\cos\left(\tfrac{2t-1}{t}r_0\right)-\cos\left(\tfrac{r_0}{t}\right)\right]}=-1.
\end{align*}
That is, $\mathsf{d}(x_0,x_1)=\pi$.
\end{proof}

\begin{theorem}\label{poles}
Assume that the metric measure space $(\mathsf{M},\mathsf{d},\mathsf{m})$ satisfies the curvature-dimension condition $\mathsf{CD}(N-1,N)$ for some  $N\ge1$ and that $\mathsf{diam}(M)\le\pi$ (which follows from the previous condition if $N\not=1$). Let $\nu$ be any optimal path measure on the spherical cone $(\Sigma(\M),\mathsf{d_\Sigma})$ satisfying $(e_0)_*\nu\ll m_N$. Then
$\nu$ gives no mass to geodesics through the poles:
$$\nu\left(\Gamma_{\mathcal S}\right)=\nu\left(\Gamma_{\mathcal N}\right)=0$$
where
$\Gamma_{\mathcal S}:=\{\gamma\in\Gamma(\Sigma(\M)): \  \gamma_t\in\mathcal S\ \mbox{\rm for some }t\in(0,1)\}$ and analogously $\Gamma_{\mathcal N}$ with ${\mathcal N}$ in the place of ${\mathcal S}$.
\end{theorem}

\begin{proof}
We follow the argumentation in the proof of assertion (iii) of Theorem 2.3. Assume that $\nu(\Gamma_{\mathcal S})>0$. Then without restriction we even may assume that $\nu(\Gamma_{\mathcal S})=1$. According to Lemma 3.4 below, for each $r\in (0,\pi)$ there exists at most one point $f(r)\in M$ such that $(f(r),r)\in\Sigma(M)$ is the initial point $\gamma_0$ of some geodesic $\gamma\in\mathsf{supp}[\nu]$.
Hence, $\mu_0:=(e_0)_*\nu$ is concentrated on the set
 $C_f:=\{(f(r),r)\in \Sigma(\M):\, r\in(0,\pi)\}$.

The curvature-dimension condition for $(\M,\mathsf d,\m)$ implies that $\m$ has no atoms and thus
$$\widehat\m_N(C_f)=0$$
which contradicts the assumption $\mu_0\ll \widehat\m_N$. Hence, $\nu(\Gamma_{\mathcal S})=0$. Analogously, we deduce $\nu(\Gamma_{\mathcal N})=0$.
\end{proof}

\begin{lemma} Under the assumptions of the previous theorem, for every $r\in(0,\pi)$ there exists at most one $x\in M$ such that $\gamma_0=(x,r)\in\Sigma(M)$ is the initial point of some geodesic $\gamma\in \mathsf{supp}[\nu]\cap \Gamma_{\mathcal S}$.
\end{lemma}

\begin{proof}
Assume $\gamma,\gamma'\in \mathsf{supp}[\nu]\cap \Gamma_{\mathcal S}$ with $\gamma_0=(x_0,r)$ and $\gamma'_0=(x'_0,r)$. According to Lemma 3.2, the fact that $\gamma$ passes through the south pole implies that $\gamma_1=(x_1,r_1)$ with $x_1\in M$ being an antipode of $x_0$, i.e. $\mathsf{d}(x_0,x_1)=\pi$.
Similarly,  $\gamma'_1=(x'_1,r'_1)$ with  $\mathsf{d}(x'_0,x'_1)=\pi$. The radii $r_1, r'_1$ are arbitrary  numbers in $(0,\pi)$.

By the very definition of $\mathsf d_\Sigma$, taking into account that the diameter of $M$ is bounded by $\pi$,
\begin{align*}
\mathsf{d}^2_\Sigma&(\gamma_0,\gamma'_1)+\mathsf{d}^2_\Sigma(\gamma'_0,\gamma_1)\\
&=\arccos^2\left[\cos r\cdot\cos r_1'+\sin r\cdot\sin r_1'\cdot\cos\mathsf{d}(x_0,x'_1)\right] \\
&\hspace{3cm}+ \arccos^2\left[\cos r\cdot\cos r_1+\sin r\cdot\sin r_1\cdot\cos\mathsf{d}(x'_0,x_1)\right]\\
&\stackrel{(\ast)}\le\arccos^2\left[\cos r\cdot\cos r_1'-\sin r\cdot\sin r_1'\right] +
 \arccos^2\left[\cos r\cdot\cos r_1-\sin r\cdot\sin r_1\right]\\
&=(r+r_1')^2+(r+r_1)^2=\mathsf{d}^2_\Sigma(\gamma'_0,\gamma'_1)+\mathsf{d}^2_\Sigma(\gamma_0,\gamma_1)
\end{align*}
with equality in $(\ast)$ if and only if $\mathsf{d}(x_0,x'_1)=\mathsf{d}(x'_0,x_1)=\pi$, that is, if and only if $x_0$ and $x_1'$ are antipodes and $x_0'$ and $x_1$ are antipodes. Therefore,
 $\mathsf{d}^2_\Sigma$-cyclical monotonicity implies $x_0=x_0'$.
\end{proof}

From now on, let us again focus on weighted Riemannian spaces, that is, $M$ is a complete, $n$-dimensional manifold equipped with its Riemannian distance $\mathsf d$ and with a measure $d\m(x)=e^{-V(x)}d\mathsf{vol}_M(x)$.
A crucial fact for our argumentation is that the punctured cone $\Sigma_0:=\Sigma(\M)\setminus\{\mathcal{S},\mathcal{N}\}$ is given as a warped product $\M\,\,{}_{\sin(r)}\!\!\times(0,\pi)$ for which the Ricci curvature can be calculated explicitly.

\begin{lemma}
(i) \ The punctured spherical cone $\Sigma_0$ is an incomplete $(n+1)$-dimensional Riemannian manifold whose tangent space $T_{({x},r)}\Sigma_0$ at $({x},r)\in\Sigma_0$ with ${x}\in\M$ and $0<r<\pi$ can be parametrized as
$T_{({x},r)}\Sigma_0=T_{x}\M\oplus\R$
and whose metric tensor is given by
$\| (v,t)\|^2_{T_{({x},r)}}=\sin^2r\cdot\| v\|^2_{T_{x}}+ t^2$
for $(v,t)\in\mathsf{T}_{({x},r)}\Sigma_0$. Furthermore, we have the equality
$$\mathsf{Ric}_{({x},r)}((v,t),(v,t))=\mathsf{Ric}_{x}(v,v)+(1-n\cos^2r)\cdot\|v\|^2_{T_{x}}+n\,t^2.$$
In particular, $\mathsf{Ric}\geq n$ on $\Sigma_0$ if and only if $\mathsf{Ric}\geq n-1$ on $M$.

(ii)\ Now let us consider the punctured $N$-spherical cone over the weighted Riemannian manifold $M$. That is, given any real $N>1$ put $W({x},r)=-(N-n)\log\sin r$ and $V({x},r)=V({x})$. Then
\begin{equation}\mathsf{Ric}^{N+1,V+W}_{({x},r)}((v,t),\, (v,t))-N\| (v,t)\|^2_{T_{{x},r}}=
\mathsf{Ric}^{N,V}_{x}(v,v)-(N-1)\| v\|^2_{T_{x}}.
\end{equation}
In particular,  $\mathsf{Ric}^{N+1,V+W}\geq N$ on $\Sigma_0$ if and only if $\mathsf{Ric}^{N,V}\geq N-1$ on $\M$.
\end{lemma}

\begin{proof}
The formula for the Ricci tensor in (i) is well-known, see \cite{neill}, Cor. 7.43, or e.g. \cite{jp}.
Note that
$$\mathsf{Ric}_{({x},t)}((v,t),(v,t))-\mathsf{Ric}_{x}(v,v)=(1-n\cos^2r)\cdot\|v\|^2_{T_{x}}+n\,t^2=
n\,\|(v,t)\|^2_{T_{({x},r)}}-(n-1)\,\|v\|^2_{T_{x}}.$$

The proof of assertion (ii) follows the lines of argumentation in the previous case of Euclidean cones, -- with appropriate modifications.
 For arbitrary $V({x},r)=V({x})$ and $W({x},r)=W(r)$ as above (depending only on the radial coordinate $r\in\R$ or on the basic coordinate ${x}\in M$, resp.) we have as before
$$\left[ \nabla V \otimes \nabla W\right]_{({x},r)} ((v,t),\, (v,t))=
{\nabla V}_x(v)\cdot W'(r)\cdot t$$ for all $({x},r)\in \Sigma_0$ and all $(v,t)\in T_{({x},r)}\Sigma_0$ and
$$\nabla W_{({x},r)}(v,t)=h'(0),\qquad
\left[  \mathsf{Hess}\, W\right]_{({x},r)} ((v,t),\, (v,t))=h''(0)
$$
where now
$$h(s):=W\left(\arccos\left[ \cos(r+st)\cdot\cos(s \cdot\sin r\cdot\| v\|_{T_{x}})\right]\right).$$
Moreover,
$$
\left[  \mathsf{Hess}\, V\right]_{({x},r)} ((v,t),\, (v,t))\quad=\quad f''(0)\quad=\quad\left[  \mathsf{Hess}\, V\right]_{{x}} (v,v)-2{\nabla V}_x(v)\cdot \cot(r)\cdot t$$
where
$$f(s)=V\left(\exp_x\left(\frac v{\|v\|_{T_x}}\cdot\arctan\frac{\tan(\sin(r)\cdot s \|v\|_{T_x})}{\sin(r+st)}\right)\right).$$
For the particular choice of $W({x},r)=-(N-n)\log \sin(r)$, some lengthy calculation yields
\begin{eqnarray*}\left[  \mathsf{Hess}\, W - \frac1{N-n}\nabla W \otimes \nabla W\right]_{({x},r)} ((v,t),\, (v,t))&=&(N-n)\left[t^2-\cos^2(r)\, \| v\|^2_{T_{x}}\right]\\
&=&(N-n)\left[\|(v,t)\|^2_{T_{({x},r)}}-\,\|v\|^2_{T_{x}}\right].\end{eqnarray*}
Hence, together with the identity from (i)
\begin{eqnarray*}
\lefteqn{\mathsf{Ric}^{N+1, V+W}_{({x},r)}((v,t),\, (v,t))}\\
&=&
\mathsf{Ric}_{({x},r)}((v,t),\, (v,t))
+\left[  \mathsf{Hess}\, (V+W) - \frac1{N-n}\nabla (V+W) \otimes \nabla (V+W)\right]_{({x},r)} ((v,t),\, (v,t))\\
&=&
\mathsf{Ric}_{x}(v,v)+n\,\|(v,t)\|^2_{T_{({x},r)}}-(n-1)\| v\|^2_{T_{x}}\\
&&\quad+\left[ \mathsf{Hess}\, V - \frac1{N-n}\nabla V \otimes \nabla V\right]_{{x}} (v,v)
+(N-n)\, \left[\|(v,t)\|^2_{T_{({x},r)}}-\,\|v\|^2_{T_{x}}\right]\\
&=&
\mathsf{Ric}^{N,V}_{x}(v,v)-(N-1)\| v\|^2_{T_{x}}
+N\,\|(v,t)\|^2_{T_{({x},r)}}
.
\end{eqnarray*}

\end{proof}

\begin{theorem}\label{sphconriemman}
Given a complete $n$-dimensional manifold $M$ equipped with its Riemannian distance $\mathsf d$ and with a weighted measure $d\m(x)=e^{-V(x)}d\mathsf{vol}_\M(x)$ for some function $V:\M\to\R$.  Then for each real number $N\ge1$ the following statements are equivalent:
\begin{itemize}
\item[(i)]
The weighted Riemannian space $(\M,\mathsf{d},\m)$ has $\mathsf{diam}(\M)\le\pi$ and satisfies the  condition $\mathsf{CD}(N-1,N)$.
\item[(ii)] The $N$-spherical cone $(\Sigma(\M),\mathsf{d_\Sigma},\m_N)$  satisfies the condition $\mathsf{CD}(N,N+1)$.
\end{itemize}
\end{theorem}

\begin{proof} This is essentially the same argumentation as in the proof of Theorem 2.5, now with Lemma 3.5 instead of Lemma 2.4.

$(i) \Rightarrow (ii)$:\  Let probability measures $\mu_0$ and $\mu_1$ on $\Sigma(\M)$ be given, absolutely continuous with respect to $\widehat\m_N$. According to Theorem \ref{poles}, any optimal path measure $\nu$ with marginal distributions $(e_0)_*\nu=\mu_0$ and $(e_1)_*\nu=\mu_1$ will give no mass to geodesics through the poles. In other words, $\nu$-almost every geodesic will  stay within the punctured cone $\Sigma_0$.

According to  Lemma 3.5, assertion (i) implies that the $(N+1)$-Ricci tensor $\mathsf{Ric}^{N+1,V+W}$ on the weighted Riemannian space $\Sigma_0$ is bounded from below by $N$.
Hence, classical arguments based on Jacobi field calculus -- exactly the same as used  to deduce Lemma 1.11 -- will imply that (1.2) holds true with $K=N$ for $\nu$-a.e. geodesic $\gamma$ which remains within $\Sigma_0$. That is, $\mathsf{CD}(N,N+1)$ holds true on $\Sigma(\M)$.

$(ii) \Rightarrow (i)$:\ The curvature-dimension condition  $\mathsf{CD}(N,N+1)$ for the $N$-spherical cone $(\Sigma(\M),\mathsf{d_\Sigma},\widehat\m_N)$  implies that this condition holds \emph{locally} on the punctured cone $\Sigma_0$. For this (non-complete) weighted Riemannian manifold, however, the local curvature-dimension condition $\mathsf{CD}_{loc}(N,N+1)$ is equivalent to the bound $\mathsf{Ric}^{N+1,V+W}\ge N$ for the $(N+1)$-Ricci tensor on $\Sigma_0$, see Lemma 1.11. Due to  Lemma 3.5, this implies $\mathsf{Ric}_M^{N,V}\ge N-1$. For the (complete) weighted Riemannian space  $(\M,\mathsf{d},\m)$, the latter in turn is equivalent to
$\mathsf{CD}(N-1,N)$.

Finally, in the case $N=1$ it remains to prove that (ii) implies the diameter bound $\mathsf{diam}(\M)\le\pi$.
This can be achieved by means of a straightforward adaptation of the argument from the proof of Corollary 2.6.
\end{proof}

\begin{corollary}
The $n$-spherical cone $(\Sigma(\M),\mathsf{d}_\Sigma,\nu)$ over a  complete $n$-dimensional Riemannian manifold $(\mathsf{M},\mathsf{d}, \V)$ satisfies $\mathsf{CD}(n,n+1)$ if and only if  $\mathsf{Ric}\geq n-1$ on $M$ and $\mathsf{diam}(\M)\le\pi$.
\end{corollary}

Theorem \ref{sphconriemman} allows to apply the Lichnerowicz theorem \cite{lva} in order to obtain a lower bound on the spectral gap of the Laplacian on the spherical cone:

\begin{corollary}[Lichnerowicz estimate, Poincar\'e inequality]
Let $(\Sigma(\M),\mathsf{d}_\Sigma,\widehat\m_n)$ be the $n$-spherical cone of a  compact $n$-dimensional Riemannian manifold $(\mathsf{M},\mathsf{d}, \V)$ with $\mathsf{Ric}\geq n-1$ and $\mathsf{diam}(\M)\le\pi$. Then for every $f\in\mathsf{Lip}(\Sigma(\M))$ fulfilling $\int_{\Sigma(\M)} f \ d\widehat\m_n=0$ the following inequality holds true,
\begin{equation*}
\int_{\Sigma(\M)} f^2d\widehat\m_n\leq \tfrac{1}{n+1}\int_{\Sigma(\M)}|\nabla f|^2d\widehat\m_n.
\end{equation*}
\end{corollary}

The Lichnerowicz estimate implies that the Laplacian $\Delta$ on the spherical cone $(\Sigma(\M),\mathsf{d}_\Sigma,\widehat\m_n)$ defined by the identity
$$\int_{\Sigma(\M)} f\cdot\Delta g \ d\widehat\m_n=-\int_{\Sigma(\M)}\nabla f\cdot\nabla g \ d\widehat\m_n$$ admits a spectral gap $\lambda_1$ of size at least $n+1$,
$$\lambda_1\geq n+1.$$

\smallskip

An analogous statement -- with $N$ in the place of $n$ -- holds true for the Laplacian on the $N$-spherical cone over a weighted $n$-dimensional Riemannian manifold satisfying $\mathsf{Ric}^{N,V}\geq N-1$.

\bigskip

\subsection*{Extension to $(\kappa,N)$-Cones}

Let us finally mention that there is a canonical extension of the concept of cones which covers both, the Euclidean cones and the spherical cones.

\begin{definition}
Given a metric measure space $(\mathsf{M},\mathsf{d},\mathsf{m})$ and numbers $\kappa\in\R, N\in (0,\infty)$ we define the $(\kappa,N)$-cone over $(\mathsf{M},\mathsf{d},\mathsf{m})$ to be the  metric measure space
$(\overline\M,\overline{\mathsf d},\overline\m)$ with

\begin{itemize}
\item[$\diamond$] $\overline\M:=\mathsf{M}\times[0,\infty)$ if $\kappa\le0$ and
$\overline\M:=\mathsf{M}\times[0,\pi/\sqrt\kappa]$ if $\kappa>0$ where all the points $(x,0)$, $x\in M$, have to be identified as well as -- in the case $\kappa>0$ -- all the points $(x,\pi/\sqrt\kappa)$.
\item[$\diamond$] For $(x,s),(y,t)\in \overline\M$
\begin{equation}\label{dk}\overline{\mathsf d}((x,s),(y,t)):=\ck^{-1}\left(\ck(s)\cdot\ck(t)+\kappa\cdot\sk(s)\cdot\sk(t)\cdot\cos\left(\mathsf{d}(x,y)\wedge\pi\right)\right)\end{equation}
where $\ck (r)=\cos(\sqrt\kappa\,r)$, $\sk (r)=\frac1{\sqrt{\kappa}}\sin(\sqrt\kappa\,r)$ if $\kappa>0$ and
$\ck (r)=\cosh(\sqrt{-\kappa}\,r)$, $\sk(r)=\frac1{\sqrt{-\kappa}}\sinh(\sqrt{-\kappa}\,r)$ if $\kappa>0$. In the case $\kappa=0$, the metric $\overline{\mathsf d}$ will be defined as in Definition 2.1. Indeed, the formula (3.2) leads in the limit $\kappa\to0$ to the definition of $\mathsf{d_{Con}}$.
\item[$\diamond$] $d\overline\m(x,s):=d\m(x)\otimes(\sk(s)^N ds)$.
\end{itemize}
\end{definition}

The metric space $(\overline\M,\overline{\mathsf d})$ obtained as such a cone over a metric space $(\M,\mathsf d)$ is discussed in detail in \cite{bi}.
In the case $\kappa=0$ it is simply the Euclidean cone and in the case $\kappa=1$ it is the spherical cone. In the case $\kappa=-1$, the cone is also called \emph{hyperbolic cone} based on $(\M,\mathsf d)$. Without too much effort, our previous results extend to the general case of $(\kappa,N)$-cones over weighted Riemannian spaces $(\mathsf{M},\mathsf{d},\mathsf{m})$. Indeed, the case $\kappa>0$ is just a rescaling of the case $\kappa=1$. Replacing all $\sin$ and $\cos$ by $\sinh$ and $\cosh$ (e.g. in Lemma 3.5) allows to switch from the case $\kappa>0$ to the case $\kappa<0$.

\begin{theorem}
Given a complete $n$-dimensional manifold $M$ equipped with its Riemannian distance $\mathsf d$ and with a weighted measure $d\m(x)=e^{-V(x)}d\mathsf{vol}_\M(x)$ for some function $V:\M\to\R$.  Then for all $\kappa\in\R$ and  $N\ge1$ the following statements are equivalent:
\begin{itemize}
\item[(i)]
The weighted Riemannian space $(\M,\mathsf{d},\m)$ has $\mathsf{diam}(\M)\le\pi$ and satisfies the  condition $\mathsf{CD}(N-1,N)$.
\item[(ii)] The $(\kappa,N)$-cone $(\overline\M,\overline{\mathsf d},\overline\m)$  satisfies the condition $\mathsf{CD}(\kappa\cdot N,N+1)$.
\end{itemize}
\end{theorem}

\
\\
\textbf{Acknowledgement.} The authors are grateful to the anonymous referee for her/his comments which helped to improve and simplify the presentation. The second author also would like to thank Jeff Cheeger for stimulating discussions during a visit at Courant Institute in 2004, in particular, for posing the problem treated in Theorem \ref{Cheeger}, as well as for valuable remarks on an early draft of this paper.

\end{document}